\documentclass{notices}
\usepackage{amsfonts,amssymb,amsmath,amscd,graphicx}

\title{Rational Approximation}

\author{
Lloyd N. Trefethen
\affil{Professor of Applied Mathematics in Residence,
School of Engineering and Applied Sciences, Harvard University}
}

\begin{document}

\maketitle

Rational approximation is an old subject, but it has gone through
a transformation in recent years with the appearance of a new
algorithm that makes computing rational approximations easy.
Indeed, many problems can now be solved to good accuracy in
milliseconds on a laptop.  This is the AAA algorithm, whose name
derives from ``adaptive Antoulas-Anderson''~\cite{aaa}.  Curiously,
the computation of polynomial approximations on domains other than
intervals and disks, which is in principle a much easier problem
but still becomes numerically intractable if a good basis is not
available, has also been transformed recently by the introduction
of the so-called Vandermonde-with-Arnoldi algorithm~\cite{VA}.

The aim of this short story is not to describe these algorithms,
which has been done elsewhere, but to show how beautifully they
can reveal the power of rational functions and the mathematics
of rational approximation.  With the help of six figures all in
the same format, we will illustrate what is sometimes called the
Walsh/Gonchar theory of this subject.  Details and references can
be found in~\cite{levin,aaa,rakh,jjiam}.

The setting is a function $f$ defined on a compact, simply connected
set $K$ in the complex plane $\mathbb C$, and $\|\cdot\|$ is the
$\infty$-norm over $K$.  (Simple connectivity doesn't matter for the
algorithm, but it makes this exposition simpler.)  We are interested
in best or near-best approximations to $f$ by rational functions
$r$ of degree $n$ (i.e., $r = p/q$ with polynomials $p$ and $q$
of degree $n$), and for comparison, by polynomials of degree $n$.
The minimax errors can be written $$ E_n = \inf_p\|f-p\|, \quad
E_{n,n} = \inf_r\|f-r\|. $$

\begin{figure}
\noindent\kern 1pt\includegraphics [clip, width=3.2in]{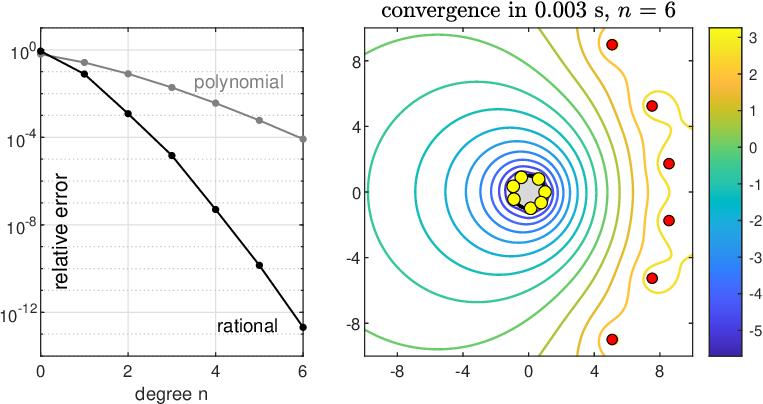}
\caption{Entire function, $f(z) = \exp(z)$}          
\end{figure}

First, consider the case where $f$ is entire.  Here it is
known that both $E_n$ and $E_{n,n}$ are guaranteed to decrease
super-exponentially as $n\to\infty$, i.e., they are both $O(C^{-n})$
for any $C>1$. Figure~1 illustrates this for approximation of $f(z)
= e^z$ on the unit disk.  On the left, we see convergence curves for
AAA rational and Vandermonde-with-Arnoldi polynomial approximation,
curving downward on this log scale.  On the right, we see a portion
of the complex plane with information about the final rational $r$
approximation to accuracy $10^{-12}$, of degree $n=6$ in this case.
The red dots mark the poles $\{\pi_k\}$ of $r$, $1\le k\le n$.  The
yellow dots mark a set of interpolation points $\{z_k\}$ with $f(z_k)
= r(z_k)$, $0\le k\le n$, that are selected by AAA in the course of
its computation.  The AAA algorithm is a ``greedy'' iteration based
on a barycentric representation of $r$ with coefficients determined
by a least-squares minimization~\cite{aaa}.  There may be other
points besides $\{z_k\}$ where $r$ also interpolates the data,
but the plots just mark these $n+1$ ``support points'' explicitly
selected step by step by AAA.  (Because of the least-squares
aspect of AAA, sample sets for $f(z)$ do not have to be determined
carefully; it is usually enough to distribute a few hundred sample
points with reasonable density along the boundary of the domain.)
The contours in this and the other plots are potential curves assoiciated with the poles
and interpolation points, as we will explain at the end.

\begin{figure}
\noindent\kern 1pt\includegraphics [clip, width=3.2in]{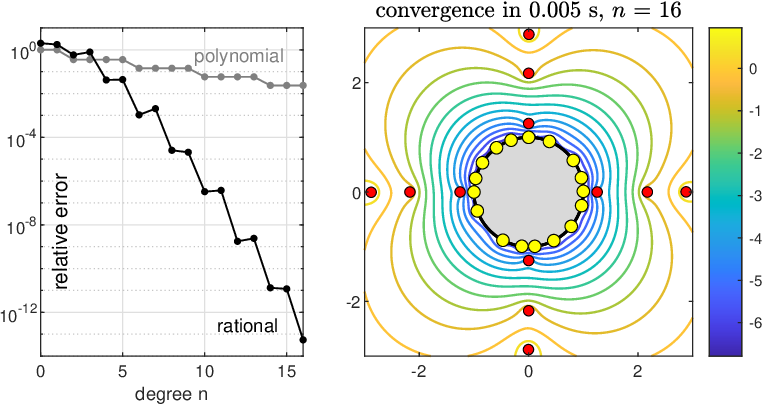}
\caption{Meromorphic function, $f(z) = \tan(z^2)$}          
\end{figure}

Next, suppose $f$ is analytic on $K$ but meromorphic in ${\mathbb
C}\backslash K$, having one or more poles.  Here rational
approximations can ``peel the poles off,'' so the decay of $E_{n,n}$
is still superexponential. Polynomials, however, are limited to
exponential convergence at a rate determined by the pole closest
to $K$ as measured by a conformal map of ${\mathbb C}\backslash K$
to the exterior of the unit disk \cite[Thm.~3]{inlets}.  Figure~2
illustrates all this for $f(z) = \tan(z^2)$ on the unit disk, with the
rational curve again curving downward but the polynomial curve now
straight.  (Four poles are off-scale, and similarly there are 4,
1, 4, and 5 poles off-scale in Figures 3, 4, 5, and 6.)

\begin{figure}
\noindent\kern 1pt\includegraphics [clip, width=3.2in]{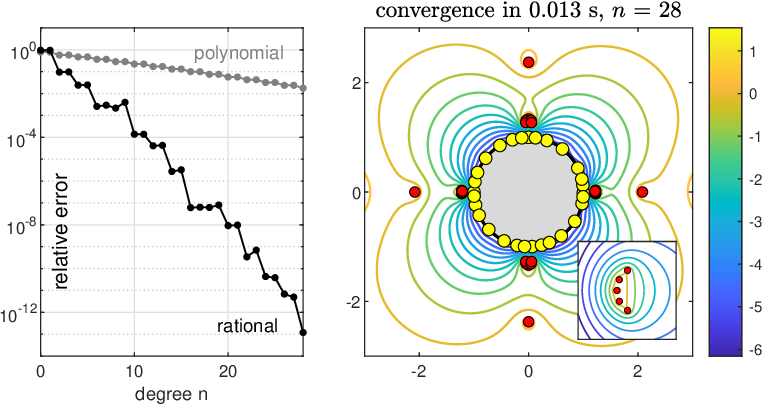}
\caption{Function with essential singularities, $f(z) = \exp(\tan(z^2))$.
Inset: zoom near $z = \sqrt{\pi/2}$}          
\end{figure}

Now suppose $f$ has essential singularities outside $K$.  You might
expect this to slow down rational approximations fundamentally,
but in fact, the convergence is still superexponential, as shown in
Fig.~3 for $f(z) = \exp(\tan(z^2))$ on the unit disk.  The contour
plot looks much as in Fig.~2, but each of the inner four red dots
is now a cluster of five poles, as illustrated in the inset.

\begin{figure}
\noindent\kern 1pt\includegraphics [clip, width=3.2in]{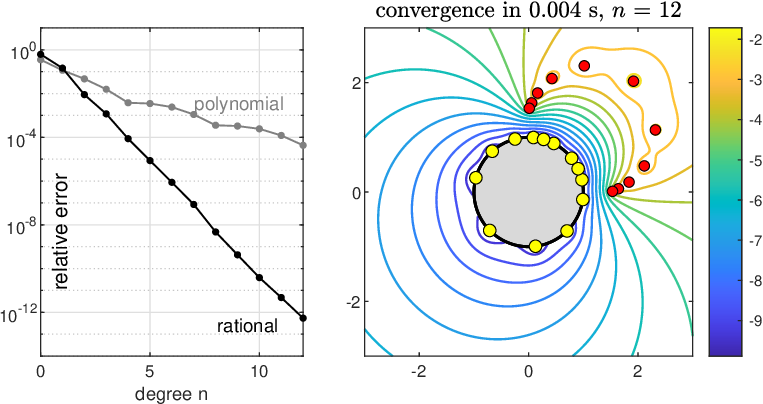}
\caption{Function with branch points, $f(z) = \sqrt{(1.5-z)(1.5\kern .4pt i-z)}$}          
\end{figure}

Branch points are a more serious obstacle to approximation, because
they are not isolated singularities. If $f$ has any branch points
outside $K$, rational approximants always slow down to exponential
convergence.  Thus both rational functions and polynomials converge
exponentially, but the rates may be very different.  I like to think
that the trouble with polynomials is that their poles are constrained
to lie at $\infty$, whereas rational functions can put poles wherever
they are most useful.  The poles of good rational approximations do
their best to approximate branch cuts, lining up beautifully along
certain curves investigated by Stahl~\cite{stahl}, as one sees in
Fig.~4 for approximation on the unit disk of a function with square
root singularities at $z=1.5$ and $1.5i$.

Finally, there are two great cases where rational functions
are far more efficient than polynomials.  The well-known one, as
discovered by Newman in 1964~\cite{newman}, is when $f$ has branch
point singularities on the boundary of $K$. Here, polynomials are
stuck with algebraic convergence, whereas rational functions can
achieve root-exponential convergence ($O(C^{-\sqrt n}\kern 1pt)$)
by clustering poles exponentially near the singularities. Figure 5
shows Newman's example of approximation of $f(z) = |z|$ on $[-1,1]$.
The rational approximations converge root-exponentially, whereas
the polynomials do no better than $O(1/n)$.

\begin{figure}
\noindent\kern 1pt\includegraphics [clip, width=3.2in]{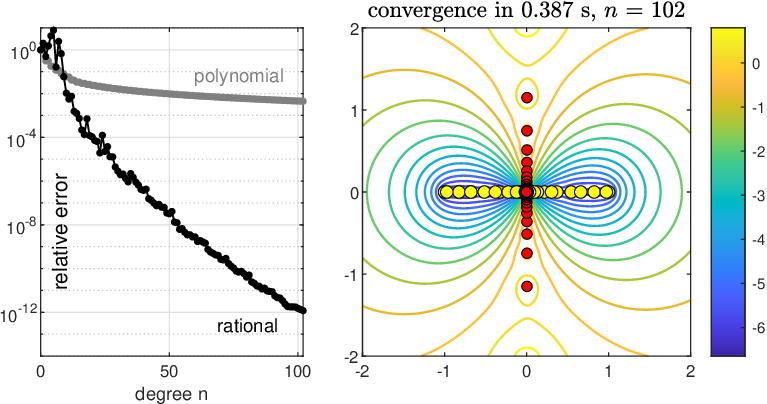}
\caption{Function with singularity on the approximation domain, $f(z) = |z|$
on $[-1,1]$}
\end{figure}

\begin{figure}
\noindent\kern -1pt\includegraphics [clip, width=3.2in]{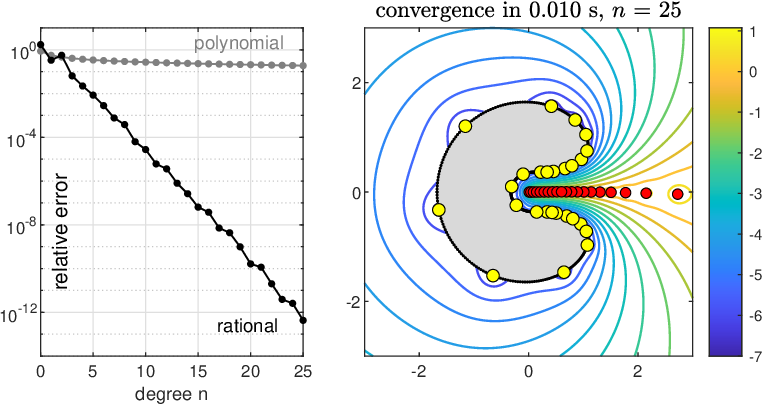}
\caption{Function analytic in a nonconvex region, $f(z) = \sqrt{-z}$}          
\end{figure}

The other case is not so well known but equally important. This
is the situation of approximation on a nonconvex domain when the
analytic continuation of $f$ has a singularity in an ``inlet'',
as will be true in most applications.  Here both polynomials and
rational functions achieve exponential convergence, $\|f-r\| =
O(C^{-n})$ for some $C>1$, but for polynomials, $C$ is
exponentially close to $1$ with respect to the length of the
inlet~\cite[sec.~7]{inlets}.  In the case shown in Fig.~6,
involving a smooth domain $K$ bent around a portion of the positive real axis,
we would need to increase the degree of $p$ asymptotically by about
$500$ for each additional digit of accuracy.  In such cases, for
practical purposes, polynomials effectively do not converge at all.
The rational functions do fine by approximating a branch cut along
$[\kern .3pt 0,\infty)$.

Now at last let us explain
the contours in the figures, which show off the Walsh/Gonchar theory that leads to
all of the results we have outlined. These contours are level
curves of a potential function $\phi$ 
which we can interpret as generated
by a positive point charge at each pole $\pi_k$
and a negative point charge at each interpolation point $z_k$,
$$ \phi(z)
= \prod_{k=0}^n (z-z_k) \left/ \prod_{k=1}^n (z-\pi_k)\right.. \eqno (1) $$
This yields an error estimate for~$r$ by a Cauchy--Hermite contour
integral exploited by Joseph Walsh beginning in the 1930\kern
.3pt s~\cite[Thm.~8.2]{walsh}, $$ f(z) - r(z) = {1\over 2\pi i}
\int_\Gamma {\phi(z)\over \phi(t)}\kern 1pt {f(t) \over t-z} \kern
1pt dt. \eqno (2) $$ This identity holds for any contour $\Gamma$
within which $f$ and $r$ are analytic, and good bounds come
when the pole and interpolation point distributions are close to
equilibrium.  Behavior as $n\to\infty$ was investigated by Gonchar
and others beginning in the 1960\kern .3pt s, and the best bounds
come with $\Gamma$ lying ``as far from $K$ as possible,'' hugging
Stahl's optimal branch cuts.  The AAA algorithm doesn't know the
Walsh/Gonchar theory, but its poles delineate the optimal branch
cuts pretty well anyway.

Something remarkable emerges if one examines the numbers on the
colorbars in our plots.  These represent powers of 10, and what is
striking is that the range of yellow to blue is typically around 7
orders of magnitude, not $12$ or more as (2) would suggest for
12-digit accuracy.  There is a factor of 2 in play here that is
related to a certain orthogonality or approximate orthogonality
along contours.  Simply put, taking absolute values in eq.\ (2)
gives a valid upper bound, but the actual error is often much
smaller because of oscillation of the integrand, leading to
cancellation of the integral to leading order.
For many near-best approximations, $E_{n,n}$ is on the order
of the square of what one would expect from (2).  Rakhmanov
speaks of the {\em Gonchar-Stahl $\rho^2$ theorem,} which details
cases in which this squaring can be guaranteed.  See~\cite{rakh}
and~\cite[sec.~7.5]{jjiam}.

The easy computations and plots of this review would not have been possible
before about 2018.  
The implications of the AAA algorithm are very wide-ranging, making
rational functions practical for numerical computation in a way
they were not before.  Exploring the possibilities is an active
research area in numerical analysis.

\bibliographystyle{plain}
\bibliography{refs}

\end{document}